\documentclass[12pt]{amsart}
\usepackage{amsfonts,amsmath,amsthm,amscd,amssymb,latexsym,cite,verbatim,texdraw,floatflt,
caption2}
\RequirePackage{hyperref}

\usepackage[a4paper]{geometry}

\theoremstyle
{plain}

\begin{document}

\title{On a question of Dikranjan  and  Zava }

\author{Igor Protasov}

\maketitle
\vskip 5pt

{\bf Abstract.} Let $G$  be a non-discrete countable metrizable abelian topological group endowed with the coarse structure  $ \mathcal{C} $  generated by compact subsets of $G$.
 We prove that 
 $asdim (G, \mathcal{C} ) = \infty$.  
For an infinite cyclic subgroup $G$  of  the circle, this answers a question of Dikranjan  and Zava [3]. 

\vskip 10pt

\vskip 10pt


\section{Introduction}

Let $X$  be a set. A family $\mathcal{E}$ of subsets of $X\times X$ is called a {\it coarse structure } if
\vskip 7pt

\begin{itemize}
\item{}     $\bigtriangleup _{X}\in  \mathcal{E}$,
$\bigtriangleup _{X}= \{(x,x): x\in X\}$;
\vskip 5pt

\item{}  if  $E$, $E^{\prime} \in \mathcal{E}$ then $E\circ E^{\prime}\in\mathcal{E}$ and
$E^{-1}\in \mathcal{E}$,   where    $E\circ E^{\prime}=\{(x,y): \exists z((x,z) \in  E,  \   \ (z, y)\in E^{\prime})\}$,   $E^{-1}=\{(y,x): (x,y)\in E\}$;
\vskip 5pt

\item{} if $E, \    E^{\prime}\in\mathcal{E}$ and 
$H \subseteq E$  then
$E\cup E^{\prime}\in  E  $  and 
$H \in \mathcal{E}$.
\vskip 5pt

\end{itemize}
\vskip 7pt

A subset $\mathcal{E}^{\prime} \subseteq \mathcal{E}$  is called a
{\it base} for $\mathcal{E}$  if, for every $E\in \mathcal{E}$, there exists
  $E^{\prime}\in \mathcal{E}^{\prime}$  such  that
  $E\subseteq E ^{\prime}$.
For $x\in X$,  $A\subseteq  X$  and
$E\in \mathcal{E}$, we denote
$E[x]= \{y\in X: (x,y) \in E\}$,
 $E [A] = \cup_{a\in A}   \   \   E[a]$
 and say that  $E[x]$
  and $E[A]$
   are {\it balls of radius $E$
   around} $x$  and $A$.

The pair $(X,\mathcal{E})$ is called a {\it coarse space} \cite{b12}.
We note that  coarse spaces defined in terms of balls were introduced under the name balleans in [8]  independently and  simultaneously  with [12], for the history see [2].

Each subset $Y\subseteq X$  defines the {\it subspace}
$(Y, \mathcal{E}_{Y})$, where $\mathcal{E}_{Y}$  is the restriction of $\mathcal{E}$  to $Y\times Y$.
A subset $Y$  is called {\it bounded} if $Y\subseteq E[x]$  for some $x\in X$  and
$E\in \mathcal{E}$.
\vskip 7pt

A family $\mathcal{F}$  of subsets  of $X$  is called $E$-{\it bounded ($E$-disjoint)} if, for each $A\in  \mathcal{F}$, there exists $x\in X$ such that $A\subseteq  E[x]$   $(E[A]\cap   B = \emptyset$  for all distinct $A, B \in \mathcal{F})$.

By the definition [12, Chapter 9],  $ asdim(X, \mathcal{E})\leq  n$  if, for each $E\in  \mathcal{E}$,   there exist $F\in \mathcal{E}$  and $F$-bounded covering $\mathcal{M}$ of $X$  which can be partitioned $\mathcal{M}= \mathcal{M} _{0}\cup  \ldots \cup \mathcal{M}_{n}$
 so that each family $\mathcal{M}_{i}$ is $E$-disjoint.
If there exists the minimal $n$  with this property then
 $ asdim(X, \mathcal{E})= n$, otherwise  $ asdim(X, \mathcal{E}) = \infty $.

Given two coarse spaces  $(X, \mathcal{E})$,  $(X^{\prime}, \mathcal{E}^{\prime})$, a mapping
   $f: X\longrightarrow  X^{\prime}$   is called {\it macro-uniform} 
(or  bornologous  \cite{b12})
 if, for each
    $E \in \mathcal{E}$,  there  exists
    $E^{\prime} \in \mathcal{E}^{\prime}$
     such that
     $f(E[x]) \subseteq  E^{\prime}[f(x)]$ for each $x\in X$.
If $f$  is a bijection such that $f$  and $f ^{-1}$  are macro-uniform  then $f$  is called an {\it asymorphism}.

Now let $G$  be a group. A family $\mathcal{I}$  of  subsets of $G$  is
 called a {\it group ideal} 
\cite{b9}, \cite{b11} 
 if $G$  contains the family
  $[G] ^{< \omega}$  of all finite subsets of $G$  and
  $A, B \in \mathcal{I}$,
   $ \  \   C\subseteq A$
    imply
   $AB ^{-1}  \in \mathcal{I}$,  $C\in \mathcal{I}$.
Every group ideal $\mathcal{I}$  defines a coarse structure on $G $
  with the base 
$\{\{ (x, y): x\in Ay \} \cup \triangle_G : A \in \mathcal{I}\}$.
We denote $G$  endowed with this coarse structure by  $(G, \mathcal{I})$.

If $G$  is discrete  then the coarse space
 $(G, [G] ^{< \omega })$  is the main subject of {\it Geometric Group Theory}, see \cite{b4}.
For coarse structures on $G$  defined by the ideal $[G] ^{< \kappa} $, where $\kappa$  is a cardinal, see \cite{b10}.

Every topological group $G$  can be endowed with a coarse structure defined by the ideal of all totally bounded subsets of  $G$.
These coarse structures were introduced and studied in \cite{b5}.
For asymptotic dimensions of locally compact abelian groups endowed with coarse structures defined by ideals of  precompact subsets see \cite{b6}.
For the coarse structure  on a topological group $G$  defined by the group ideal  generated by converging sequences, see 
 \cite{b7}.

For a topological group
$G$,
 we denote by  $\mathcal{C}$  the group ideal of precompact  subsets of $G$ ($A$ is precompact if $cl A$ is compact).
In  [3, Problem 5.1]  Dikranjan and Zava observed that $asdim (G,\mathcal{C} )> 0$   for an infinite cyclic subgroup of the circle and asked about the value of  $asdim (G,\mathcal{C} )$.

\section{Results }

We denote by $\Phi$ the family of  all  mappings $\phi: \omega \longrightarrow [\omega]^{<\omega} $
 such  that, for each  
$n\in \omega $, $n\in \phi (n)$  and $    \{ m: n\in \phi (m)\}$  is finite. 
We consider the family 
$\mathcal{F}$ of all subsets of $\omega \times \omega$ of the form
$$\{ (n, k):  k\in  \phi (n), \  n\in \omega  \},  \ \phi\in  \Phi $$
and note that 
 $\mathcal{F}$  is a coarse  structure  on $\omega $, see [11, Example 1.4.6]. 
The universal property of $(\omega ,\mathcal{F}): $
 if $\mathcal{E}$ 
 is a coarse structure on $\omega$ such that every bounded subset in
 $(\omega,\mathcal{ E})$ is finite then each injective mapping 
$f: (\omega, \mathcal{E})  \longrightarrow  (\omega, \mathcal{F}) $ is 
macro-uniform.

\vspace{7 mm}

{\bf Theorem 1}. {\it $Asdim   (\omega, \mathcal{F})  = 1$.

\vspace{5 mm}

Proof.}  We take an arbitrary 
$F\in \mathcal{F}$
 such that   $F= F^{-1}$, $(n, n+1)\in F$, $n\in\omega$. 
We put $P_0 = F[0]$,  
 $P_{n+1}  = F[P_n]\setminus (P_0 \cup \dots \cup P_n)$
and note that 
$\cup_{n<\omega} P_n = \omega $.  We denote
 $\mathcal{A}_0 = \{ P_{2n}: n< \omega \}$, 
$\mathcal{A}_1 = \{ P_{2n+1}: n< \omega \}$
and  observe that $\mathcal{A}_0$,  $\mathcal{A}_1$ are $F$-disjoint.

We define a mapping  $\phi\in \Phi$ by $m\in \phi (n)$  if and only if $m\in P_n$. 
We put $H=\{(n,m): m \in  \phi(n)$, $n<\omega \}$ and note that the family $\{ P_n : n< \omega \}$ is $H$-bounded so $asdim (\omega, \mathcal{F})\leq 1.$
Since $(n, n+1)\in F$,  $\omega$ can not be partitioned into $F$-disjoint uniformly bounded subsets, so $asdim (\omega,  \mathcal{F} ) \neq 0$ and
$asdim (\omega, \mathcal{F})= 1.$
 $ \  \Box$

\vspace{7 mm}

{\bf Theorem 2}. {\it   Let $G$ be a subgroup of a topological group $H$
 such that there exists an injective sequence  $(a_n)_{n<\omega}$  in $G$  converging  to $h\in H\setminus G$, $A=\{ a_n : n<\omega\}$. 
Then the subspace $A$ of 
$ (G,\mathcal{C} )$
 is asymorphic  to 
$(\omega,  \mathcal{F} ) $.

\vspace{5 mm}

Proof.}
We show that the mapping
$f: (\omega, \mathcal{F})  \longrightarrow  A, $ 
$f(n)=a_n$
is an asymorphism.
Since $h\notin G$, each bounded subset of $A$ is finite. 
By the universal property of
$(\omega,  \mathcal{F} ) $,  $f^{-1}$ 
is macro-uniform. 

To prove that $f$ is macro-uniform, we  take an arbitrary mapping $\phi: A\longrightarrow [A]^{<\omega} $ such that $a_n \in \phi (a_n)$  and
 $\{a_ m: a_ n \in \phi (a_m)\}$
 is finite. 
We put $K=\cup _{n<\omega}  \phi (a_n) a_n ^{-1} $.
By the choice of $\phi$, every injective sequence in $K$ converges to the identity of $G$.  
Hence, $K$ is compact and $\phi (a_n)\subseteq K a_n$.  
$ \  \Box$

\vspace{7 mm}

{\bf Theorem 3}. {\it  For every
non-discrete
 countable 
metrizable abelian topological group  $G$,  $asdim (G, \mathcal{C}) = \infty $.

\vspace{5 mm}

Proof.} 
We fix a natural  number $m$  and  prove that  
$asdim (G, \mathcal{C}) \geq m$.
Since the completion  $H$  of  $G$  is uncountable, we can choose  $h_1, \dots ,  h_m \in H\setminus G$  such that 
\vspace{6 mm}

$(1)   \  \   \    i_1 h_1  + \dots +  i_m  h_m  \in  G,  ( i_1, \dots , i_m ) \in \{ -1, 0, 1  \}^m $ implies $ i_1 = \dots = i_m=0.$

\vspace{6 mm}

Then we choose injective sequences  
 $(a_{1n} ) _{n<\omega}, \dots ,      (a_{mn} ) _{n<\omega} $
converging to $h_1, \dots ,  h_m $ such that 

\vspace{6 mm}

$(2)   \  \   \   a_{1 i_1 } + \dots + a_{m  i_m}=   a_{1 j_1 } + \dots + a_{m j_m}$
 implies $ i_1 =j_1,  \dots,   i_m=j_m.$

\vspace{6 mm}

By  $(2), $ 
for 
$A= A_1   + \dots + A_m$,  $ A_i=\{ a_{in} : n<\omega \}$,
 the mapping  
$f: A\longrightarrow (\omega,  \mathcal{F})^m$, 
 $f( a_{1 i_1 } + \dots + a_{m i_m}= (i_1, \dots, i_m)$
is bijection. 
By Theorem 1 and [1],  
$asdim (\omega,\mathcal{F} )^m \leq m$,    so
  it suffices to  prove that $f$  is an asymorphism.
Applying Theorem 2, we see that $f^{-1}$  is  macro-uniform. 

To prove that $f$  is macro-uniform,  we  take  an arbitrary  compact  $K$ in $G$,  $0\in K$ and show that  there exist
$\phi_1:  A_1\longrightarrow  A_1 ^{<\omega}, \dots  ,  \phi _m :  A_m\longrightarrow  A_m ^{<\omega}$
  such that

\vspace{6 mm}

$(3)   \  \   \   a\in  \phi_{ i }(a), $ and $\{ b\in A_i : a\in \phi _i (b)\}$
is  finite,  $a\in A_i$, $i\in \{1, \dots, m\}$;

\vspace{6 mm}

$(4)   \  \   \    $   $ A\cap (K +  a_{1 i_1} + \dots + a_{m  i_m})\subseteq   \phi_{1 } (a_{1 i_1} ) + \dots +  \phi_{m }  (a_{m i_m}).$

\vspace{6 mm}

We fix $k\in \{1,  \dots, m \}$, $j\in \omega$  and denote by $\phi_k (a_{kj})$  the set of all $a_{ks}\in A_k$  such that 
$( A_1   + \dots + A_{k-1}+ a_{kj} + A_{k+1} + \dots + A_m) \cap
 (K+ A_1   + \dots + A_{k-1}+ a_{ks} + A_{k+1} + \dots + A_m )  \neq\emptyset.$

Applying $(1)$, we conclude that 
 $\phi_k  ( a_{kj})$
  is finite and $\phi_1 , \dots, \phi_m$ satisfy $(3)$,   $(4)$.
$ \  \Box$

\vspace{6 mm}

CONTACT INFORMATION

I.~Protasov: \\
Faculty of Computer Science and Cybernetics  \\
        Kyiv University  \\
         Academic Glushkov pr. 4d  \\
         03680 Kyiv, Ukraine \\ i.v.protasov@gmail.com

\end{document}